\documentclass[a4paper,12pt]{amsart}%{article}
\usepackage{framed}

\usepackage[T1]{fontenc}
\usepackage{amsmath,amssymb,amscd,amsfonts,amsthm,mathrsfs,eucal,dsfont,verbatim}
\usepackage{graphicx}
\usepackage{subcaption}

\usepackage{lscape,enumerate}
\usepackage{color,soul}

\usepackage[margin=2in]{geometry}
\usepackage{hyperref}
\usepackage{marginnote}%this avoids floating marginal notes

\DeclareMathOperator{\Vol}{Vol}

\newcommand{\real}{\mathds{R}}

\newcommand{\tr}{\operatorname{tr}}

\newcommand{\Ee}{\mathds E}
\newcommand{\Pp}{\mathds P}
\newcommand{\I}{\mathds 1}

\newcommand{\CC}{\mathcal{C}}
\newcommand{\LL}{\mathcal{L}}
\newcommand{\Aa}{\mathcal{A}}

\def\1{1\!\!\hbox{{\rm I}}}

\newcommand{\HH}{\mathcal{H}}

\newcommand{\rd}{{\mathds{R}^d}}

\newcommand{\dist}{\operatorname{dist}}

\usepackage{color}%for comments

\evensidemargin
-1cm

\usepackage[dvipsnames]{xcolor}
\usepackage{xspace}
\usepackage{url}  
\usepackage{listings}  
\usepackage{fancyvrb}
\usepackage{multirow}
\usepackage{array}
\usepackage[T2A]{fontenc}
\usepackage[english]{babel}

\numberwithin{equation}{section}
\theoremstyle{plain}
\newtheorem{theorem}{Theorem}

\newtheorem{proposition}[theorem]{Proposition}

\theoremstyle{definition}
\newtheorem{remark}[theorem]{Remark}

\newtheorem{example}[theorem]{Example}

\numberwithin{theorem}{section}
\numberwithin{definition}{section}
\begin{document}                 
	
		\title{\bfseries Weak error of approximation of  a L\'evy-type process using dynamic cutting}
	
	\author{V.\ Knopova}

	\address{Institute of Mathematics, NAS of Ukraine, Tereschenkivska st.~3, 01024,  Kyiv, Ukraine.  \normalfont\url{vicknopova@gmail.com}}
	
	\thanks{The authors thanks  VUIAS for the financial support}
	
	\maketitle 
	
	\begin{abstract} \noindent
		In this note we consider, on one hand, the  implementation of  the   parametrix method for  estimation of the distance between functionals  $\Ee^x f(X_t)$ and $\Ee^x f(\widetilde X_t)$ of two L\'evy-type processes $X$ and $\tilde X$. On the other hand, we estimate the weak error appearing by substituting small jumps by a (time-inhomogeneous) diffusion. Both methods rely on the time-dependent cutting procedure, designed to remove small jumps in a time-dependent way.
		
		\medskip\noindent
		\textbf{Key words:}
		Weak error; L\'evy-type processes,  parametrix construction, dynamic cutting.
		
		\smallskip\noindent
		\textbf{MSC2020:}
		65C30, 60H10; 35R09; 60H35;  60H30;  60G52.
	\end{abstract}             %%
%%%%%%%%%%%%%%%%%%%%%%%%%%%%%%%%%%%%%%%%%%%%%%%%%%%%%%%%%%%%%%%%%%%

	\section{Introduction}
	
	Weak sense approximations of solutions to SDEs driven by L\'evy processes become a popular topic during the last decades.  The problem in simulating of the solution and estimating its functionals comes from the fact, that in general   L\'evy processes exhibit infinite activity. Although the small jumps can be cut away from the trajectory of the driving L\'evy process, the other possible way might be to substitute the small -jump part by an independent Wiener process, which seem to be first proposed by Asmussen and Rosi\'nski in \cite{AR01}.  Since  then the question of efficiency of this  approximation is risen in many papers, see, in particular,  \cite{DMT21}, \cite{BM25}.

	In this note we investigate how the time-dependent cutting of the L\'evy kernel can be applied to the approximation of  functionals of type $\Ee^x f(X_t)$  for suitable $f$.  The key idea used in the paper is a special time-dependent cutting of the respective L\'evy(type) measures.  
	The idea of a time-dependent cutting originated in earlier works  \cite{KK13}, \cite{KK18}, \cite{KKS21},  where the small-time estimates for a L\'evy process where constructed and later on used for the paramtrix construction of the transition probability density of a L\'evy-type process. The usage of the time-dependent or, as one might  call it,  ``dynamic'' cutting, is inevitable in such a construction -- the ``dynamic'' way allows to chose the correct  zero-order approximation (denoted below by $p^0$, see \eqref{p01}) and to show the convergence of the series, representing the solution (cf. \eqref{wt1}), whereas  the analogous construction with the constant-level cutting may fail to produce such a result. The  reasonable question is whether this construction may appear to be useful for estimating the weak error in estimation of functionals of type $\Ee^x f(X_t)$, where $X$ is the respective L\'evy(-type) process. In Theorem~\ref{T1} we show, that if one can constrict by the parametrix method with the \emph{same} zero order approximation the transition probability densities  of  processes  $X$ and $\tilde X$, respectively, then under certain assumptions on the series one gets  the upper bound on the weak error $|\Ee^x f(X_t) - \Ee^x f(\tilde X_t)|$. In the applied parametrix method the dynamic cutting plays the crucial role.  Further, we generalize this result in Theorem~\ref{t3}, where under some symmetry assumption and absence of the diffusion part, the weak error can be bounded from above by the H\"older-type distance between the respective normalized kernels.
	
	Our next results concern the dymanic cutting of small jumps in the weak error approximation, without applying the parametrix method.   While the usual constant level cutting  has  such a definite advantage as its simplicity, we observed that cutting jumps in a time-dependent way also  may have some advantages, in particular, if the function $f$ in $\Ee^x f(X_t)$  lacks regularity. In Theorem~\ref{t2} we provided the upper bound on the weak error under the assumption that there exists a classical solution to the Cauchy problem, involving the L\'evy-type operator, related to $X$ via the martingale problem. We pointed out in Remark~\ref{rem2}, that the dynamic cutting allows to take $f$ from the H\"older class $\CC^\eta$ with smaller  $\eta$ as that allowed by the constant level cutting.
	
	We would like to note that  the papers  \cite{KK14}, \cite{GHKK15}, \cite{GHK18}  are close to the presented results in the sense that they contain the results on the Wasserstein distance between the paths of the solutions of L\'evy-driven SDEs,  but our motivation and the used methodology is  different. 
	On the other hand, there are results on application of the  parametrix method  for investigating  stability of transition densities  under perturbations of  the coefficients of diffusion processes, see \cite{BK21} and \cite{KKM17}, and the approximatin error, see  \cite{KM11}, \cite{Hu16}. However,  the parametrix constriction used  in \cite{KKM17} contains two different convolution schemes, and is a bit apart from the framework of Theorem~\ref{T1}. 
	
	In order to keep the presentation as transparent as possible, we formulated the results  in a rather simple situation, e.g.  the L\'evy(-type) kernel being in a sense comparable to a symmetric $\alpha$-stable one, absent drift (in Theorems~\ref{T1} and \ref{t3}) and the case of bounded coefficients. These results allow extension, which we plan to accomplish in the subsequent papers.

	\section{Settings and results}

	\subsection{Parametrix construction} 
	Let $\LL$ be an integro-differential operator of the form 
	\begin{equation}\label{LL}
		\LL f (x) := b(x)\cdot D f(x) + \int_\rd \left(f(x+u) - f(x) - D f(x) \cdot u \I_{|u|\leq 1}\right)\nu^x(du), 
	\end{equation}
	where the function  $f$ is from the space $C_0^2(\rd)$ of twice continuously differentiable functions, having compact support, $b: \rd\mapsto \rd$ is some bounded  measurable function, and  $\nu^x(\cdot)$ is a L\'evy-type kernel, i.e. 
	$ \int_\rd (1\wedge |u|^2)\nu^x(du)<\infty$ for any $x\in \rd$. 
	We specify the assumptions on $b$ and $\nu^x(du)$ below.  Our standing reference for operators of type \eqref{LL} is \cite{Ja01}. 
	Through the paper we assume that

	\begin{center}
		\begin{itemize}
			\item[\textbf{(M)}]  The martingale problem for $(\LL, C_0^\infty(\rd))$ is well-posed in the space of c\'adl\'ag functions $D([0,T],\rd)$.
		\end{itemize}
	\end{center}
	
	Recall the  parametrix construction (see, e.g., \cite{Fr64}). Let $p=p_{s,t}(x,y)$ be the fundamental solution to 
	\begin{equation}\label{fund1}
		(\partial_s + \LL_x) p =0, 
	\end{equation}
	that is, $p$ solves \eqref{fund1} and  for any $f\in C_\infty(\rd)$  one has
	$\lim_{s\uparrow t} \int_\rd f(y) p_{s,t}(x,y)dy =f(x)$. 
	Here the notation $\LL_x$ means that the operator $\LL$ is applied w.r.t. variable $x$; we drop the subscript $x$ when it does not cause misunderstanding.   Write $p$ as 
	\begin{equation}\label{p01}
		p= p^0+ r, 
	\end{equation} 
	where $p^0 = p^0_{s,t}(x,y)$ is ``some'' the zero-order approximation and $r=r_{s,t}(x,y)$ is the remainder.  Suppose that $p^0$ solves the equation  $
	\left(\partial_s+\LL^{0,(y)}_x\right)p^0=0$
	for some operator  $\LL^{0,(y)}$, whose coefficients are frozen an the end-point $y$.  Under the initial  heuristic assumption that $p$  solves \eqref{fund1} we derive 
	\begin{align*}
		0 &= (\partial_s + \LL_x) p = (\partial_s + \LL^{0,(y)}_x)p^0 + (\LL_x- \LL^{0,(y)}_x)p^0 + (\partial_s + \LL_x)  r =(\LL_x- \LL^{0,(y)}_x)p^0  + (\partial_s + \LL_x)  r, 
	\end{align*}
	that is, the remainder $r$ is the solution to 
	$(\partial_s + \LL_x)  r= - (\LL_x- \LL^{0,(y)}_x)p^0$. 
	Denote  now 
	\begin{equation}\label{H}
		H:= (\LL^{0,(y)}_x- \LL_x)p^0. 
	\end{equation}
	Then using again the heuristic assumption that $p$ is the fundamental solution, we derive $r= p\otimes H$, 
	where $\otimes$ denotes the space-time convolution. 
	Iterating this identity  yields 
	\begin{align*}
		p = p^0 + (p^0 + p \otimes H)\otimes H = p^0 + p^0 \otimes H + p\otimes H^{\otimes 2} = \dots = p^0 \otimes \left(\sum_{k=0}^\infty H^{\otimes k} \right), 
	\end{align*}
	or 
	\begin{equation}\label{wt1}
		p = p^0 + p^0 \otimes \HH, 
	\end{equation}
	where  $\HH := \sum_{k=1}^\infty H^{\otimes k}$. 
	The respective zero-order approximation $p^0$ is called the \emph{parametrix}, and \eqref{wt1} is often called the \emph{parametrix representation} of $p$.

	Similarly,  let $\widetilde \LL$ is of the form \eqref{LL}, and let $\tilde{p}$ be the fundamental solution to 
	\begin{equation}\label{fund11}
		(\partial_s + \widetilde \LL_x) \tilde p =0, 
	\end{equation}
	and that we can write the parametrix representation for $\tilde p$ 
	\begin{equation}\label{wt2} 
		\tilde 	p = p^0 + p^0 \otimes \widetilde \HH, 
	\end{equation}
	with the same parametrix $p^0$, where $ \widetilde	\HH := \sum_{k=1}^\infty \widetilde H^{\otimes k}$, 
	and $ 	\widetilde H:= (\LL^{0,(y)}_x- \widetilde \LL_x)p^0$. 
	Then one has 
	\begin{equation}\label{ppdif}
		p-\tilde p = \sum_{n=1}^\infty p^0 \otimes \left( H^{\otimes n}- \widetilde H^{\otimes n}\right).
	\end{equation}
	Note, that  for  every $n\geq 1$, the following telescoping identity holds:
	\begin{equation}\label{tele-eq10}
		H^{\otimes n}-\widetilde H^{\otimes n}     = 
		\sum_{k=0}^{n-1}
		H^{\otimes k}
		\otimes
		(H-\widetilde H)
		\otimes
		\widetilde H^{\otimes(n-1-k)}.
	\end{equation}
	Moreover, by the definitions of $H$ and $\widetilde H$,
	\begin{equation}\label{eqHdif}
		H-\widetilde H=
		(\widetilde{\LL} - \LL)p^0.
	\end{equation}
	Thus, \eqref{ppdif} can be rewritten as 
	\begin{equation}\label{diffL}
		p-\widetilde p=
		\sum_{n=1}^{\infty}
		\sum_{k=0}^{n-1}
		p^0\otimes
		H^{\otimes k}
		\otimes
		\bigl((\widetilde{\LL}- \LL)p^0\bigr)
		\otimes
		\widetilde H^{\otimes(n-1-k)}.
	\end{equation}
	This representation will play the key role in obtaining the upper bound on the  weak error.
	
	\subsection{Weak error I: application of the parametrix method}
	
	Let  $(\LL, C_0^2)$ and $(\widetilde \LL, C_0^2)$ be operators of type \eqref{LL}, and    denote by $X$ and $\tilde{X}$ the respective processes.  Denote also by $C_\infty(\rd)$ the space of continuous functions, vanishing at infinity. 
	
	Below we present the main result of this section.
	
	\begin{theorem}\label{T1}
		Suppose that $X$ and $\tilde X$  admit  transition probability densities, denoted by  $p$ and $\tilde p$, respectively,  which in turn admits the parametrix representations  \eqref{wt1} and \eqref{wt2}, with the same zero-order approximation $p^0$. 
		Suppose that the coefficients in $\widetilde \LL$ depend on a parameter $h\in (0,1)$ in such a way, that for some $\kappa, \,\delta_1, \delta_2>0$, the estimates below are satisfied: 
		\begin{equation}\label{LLeq}
			\sup_x\left\|
			(\widetilde \LL_x- \LL_x)p^0_{s,t}(x,\cdot)
			\right\|_{L_1}
			\leq
			C h
			(t-s)^{-1+\delta_1},
		\end{equation}
		
		\begin{equation}\label{HHeq}
			\sup_x\left(\|H_{s,t}(x,\cdot)\|_{L_1}
			+
			\|\widetilde H_{s,t}(x,\cdot)\|_{L_1}\right)
			\leq
			C 
			(t-s)^{-1+\delta_2}. 
		\end{equation}
		Then for $f\in C_\infty(\rd)$, $\delta:= \min(\delta_1, \delta_2)$,  we have 
		\begin{equation}
			\sup_x	\left|
			\Ee^x f(X_t)- \Ee^x f(\widetilde X_t)
			\right|
			\leq
			C_T h t^\delta\|f\|_\infty, \quad 0<t\leq T. 
		\end{equation}
	\end{theorem}
	\begin{remark}
		a) 	Note that the order of convergence does not depend on a particular zero-order approximation, provided that Assumptions of Theorem~\ref{T1} are satisfied.  b) We  require the conditions on the zero order approximation of the density, not on the density itself, which allows   more flexibility than required by using directly the Ito formula in order to get the weak error.  
	\end{remark}
	\begin{proof}[Proof of Theorem~\ref{T1}]
		In order not to overload the notation, we suppress the parameter $h$ when we write the process, related to $\widetilde \LL \equiv \widetilde \LL^h$,  and keep the same strategy for its  transition probability density.

		Since $p_{0,t}$ and $\widetilde p_{0,t}$ are the transition probability densities of $X$ and $\widetilde X$, respectively, we have 
		\begin{equation}\label{eq:weak_error_density}
			\begin{aligned}
				\mathbb{E}^x f(X_t)-
				\mathbb{E}^x f(\widetilde X_t)
				&=
				\int_{\rd}
				\bigl[
				p_{0,t}(x,y)-\widetilde p_{0,t}(x,y)
				\bigr] 	f(y)dy.
			\end{aligned}
		\end{equation}
		The key quantity in the comparison of the two parametrix expansions and getting the weak error  is \eqref{eqHdif}.  Take  $s=0$ and apply \eqref{diffL}: 
		$$
		p_{0,t}(x,y)-\widetilde p_{0,t}(x,y)=
		\sum_{n=1}^{\infty}
		\sum_{k=0}^{n-1}
		p^0\otimes
		H^{\otimes k}
		\otimes
		\bigl((\widetilde{\LL}- \LL)p^0\bigr)
		\otimes
		\widetilde H^{\otimes(n-1-k)}_{0,t}(x,y).
		$$
		Using the telescoping identity \eqref{tele-eq10}
		together with \eqref{eqHdif}, we can observe the following effect. For a fixed $n$ and $k$, the factor $(H-\widetilde H)$ contributes
		$h$, whereas the remaining $n-1$ parametrix kernels  do not contribute in $h$. Hence, using \eqref{LLeq}, \eqref{HHeq}, we obtain
		
		$$
		\left\|
		p^0
		\otimes H^{\otimes k}
		\otimes (H-\widetilde H)
		\otimes\widetilde H^{\otimes(n-1-k)}_{0,t}
		(x,y)
		\right\|_{L_1}
		\leq
		hC^n
		\frac{t^{n\delta}}
		{\Gamma(1+n\delta)}, 
		$$
		where $\delta:= \min(\delta_1, \delta_2)$. Since there are $n$ possible positions of $(H-\widetilde H)$, it follows that
		\begin{equation}
			\left\|
			p_{0,t}(x,\cdot)-\widetilde p_{0,t}(x,\cdot)
			\right\|_{L_1}
			\leq
			C h
			\sum_{n=1}^{\infty}
			\frac{n(C t^\delta)^n}
			{\Gamma(1+n\delta)}\lesssim C_T h t^\delta, \quad 0< t\leq T, 
		\end{equation}
		yielding 
		\begin{align*}
			\left|
			\Ee^x f(X_t)-\Ee^x f(\widetilde X_t)
			\right|
			&\leq
			\|f\|_\infty
			\int_{\mathbb R^d}
			\left|
			p_{0,t}(x,y)-\widetilde p_{0,t}(x,y)
			\right|\,dy
			\leq
			C_T h t^\delta\|f\|_\infty.\qedhere 
		\end{align*}
		
	\end{proof}
	
	Let us provide the sufficient conditions,  under  which the assumptions of   Theorem~\ref{T1}  hold true.  As we already mentioned in the introduction, for the sake of transparency  we do not pursue here the most general case, and assume that both operators $\LL$ and $\widetilde \LL$ from Theorem~\ref{T1} possess no drift and no quadratic parts.
	
	Take a L\'evy-type kernel $\nu^x(du)$ satisfying the following assumptions: 
	
	\medskip 
	
	\begin{itemize} 
		
		\item[\textbf{(N1)}]   $\nu^x(du)$ is symmetric in $u$ and  there exists $r_0>0$, such that  for any $x\in \rd$, $r\in (0,r_0]$, 
		\begin{equation}\label{compare}
			\nu^x(\{|u|>r\})\asymp r^{-\alpha}; 
		\end{equation}
		\item[\textbf{(N2)}]
		\begin{equation}\label{compare2}
			\inf_{\ell \in S^d} \int_{|u|\leq r} (u\cdot \ell)^2 \nu^x(du) \gtrsim  r^{2-\alpha}, 
		\end{equation} 
		where $S^d$ is a unit sphere in $\rd$.

		\item[\textbf{(N3)}] There exists $\eta\in (0,1]$, such that for any $r>0$ one has
		\begin{equation}\label{Lip}
			\int_{|u|\leq r}   |u|^2 |\nu^x(du) - \nu^y(du)|\lesssim (|x-y|^\eta \wedge 1) r^{2-\alpha}, \quad x,y\in \rd. 
		\end{equation}
	\end{itemize}
	
	\medskip 
	
	Let
	\begin{equation}\label{L1} 
		\LL f(x):= \frac12 \int_\rd 
		\bigl[
		f(x+u)+f(x-u)-2f(x)
		\bigr]\nu^x(du),\quad f\in C^2_0(\rd). 
	\end{equation}
	and take $\widetilde  \LL$ in the form
	\begin{equation}\label{Lb}
		\widetilde \LL f(x) :=  \LL  f(x)+ 
		h^\kappa B f(x), 
	\end{equation}
	\begin{equation}\label{B}
		Bf(x) := \int_\rd 
		(f(x+u)+f(x-u)-2f(x)) \frac{b(x,u)du}{|u|^{d+\beta}}, \quad f\in C^2_0(\rd),  
	\end{equation}
	where $b(x,u)$ is a positive, symmetric in $u$ bounded measurable function, and $\beta \in (0,\alpha)$. 
	\begin{proposition}\label{Pr1}
		Let $\LL$ and $\widetilde \LL $ be defined in \eqref{L1} and \eqref{Lb}, respectively, and suppose that assumptions \textbf{(N1)} -- \textbf{(N3)} are satisfied.  Then 
		
		\begin{itemize}
			
			\item[(a)] There exists $p^0$ such that the parametrix series \eqref{wt1} and \eqref{wt2} converge in $L_1$; 
			
			\item[(b)]  The martingale problem for $(\LL, C_0^2)$ (resp., $(\widetilde \LL, C_0^2)$) is well posed;

			\item[(c)] The assumptions \eqref{LLeq} and \eqref{HHeq} of Theorem~\ref{T1} are satisfied. 
		\end{itemize}
		
	\end{proposition}
	
	\begin{proof}
		\emph{(a)}   First we specify  the zero order approximation for both parametrix series.  Let 
		$$ 	\psi_t (z,\xi):= \int_{|u|\leq t^{\epsilon/\alpha}} (1-e^{i\xi u}+ i\xi u)\nu(z,du),
		$$ 
		and 
		$$p_t^{(z)}(x):= (2\pi)^{-d}\int_\rd e^{-i\xi x -\int_0^{t} \psi_r(z,\xi)dr }d\xi
		.
		$$
		Note that  \eqref{compare} implies  (see \cite{KK18})  that $\psi_t(z,\xi) \asymp |\xi|^\alpha + const$, $\xi\in \rd$, which in turn yields the existence of  $p_t^{(z)}(x)$.   Here and in what follows we omit the parameter $\epsilon$ when it does not cause misunderstanding. 
		
		One can show (see \cite{KK18}, \cite{KKS21} for details) that $p_{t}^{(z)}$ satisfies $ 	\left(\partial_s + \LL^{0, z,t-s}\right)p_{t-s}^{(z)}(x)=0$,
		where 
		\begin{equation}\label{L11}
			\LL^{0,z,t}_x f(x) := \frac12
			\int_{|u|\leq t^{\epsilon/\alpha} }
			\bigl[
			f(x+u)+ f(x-u)-2f(x)
			\bigr]\nu(z,du). 
		\end{equation}
		Similarly to as it was done in \cite{KK18}, \cite{KKS21}, one can show that 
		\begin{equation}\label{difx}
			\sup_z |D^k_x p_t^{(z)} (x)| \lesssim t^{-(d+k)/\alpha} e^{- c|x| t^{-\epsilon/\alpha}},\quad x\in \rd, \quad 0\leq t\leq T, 
		\end{equation}
		\begin{equation}\label{dift}
			\sup_z|\partial_t p_t^{(z)} (x) | \lesssim   \max(t^{-1}, t^{-1/\alpha})t^{-d/\alpha} e^{- c|x| t^{-\epsilon/\alpha}},\quad x\in \rd, \quad 0\leq t\leq T. 
		\end{equation}
		Put $	p^0_t (x,y):= p_t^{(y)} (y-x)$. 
		Following the strategy from \cite{KK18}, \cite{KKS21}, we decompose the function $H$ defined in \eqref{H} in the following way:
		
		\begin{align*}
			H_{s,t}(x,y)& = \int_{|u|\leq (t-s)^{\epsilon/\alpha}} 
			\left(\int_{-1}^1  (1-|\theta|) u' D^2 p_{t-s}^0(x+\theta u,y)u d\theta \right) (\nu^x(du) - \nu^y(du)) \\
			& \quad +\int_{|u|\geq (t-s)^{\epsilon/\alpha}} [ p^0_{t-s}(x+u,y)- p^0_{t-s}(x,y)]\nu^x(du) = I_1+I_2. 
		\end{align*}
		Using the gradient estimate  
		\eqref{difx} on $p^{(y)}$ and the Lipschitz property \eqref{Lip},  we  derive 
		\begin{align*}
			|I_1|& \lesssim  \int_{|u|\leq (t-s)^{\epsilon/\alpha}}  (t-s)^{-\frac{d+2}{\alpha}} 
			\left(\int_{-1}^1  (1-|\theta|) e^{-c|y-x-u|(t-s)^{-\epsilon/\alpha}}d\theta \right) |u|^2 |\nu^x(du) - \nu^y(du)|\\
			& \lesssim  (t-s)^{-\frac{d+2}{\alpha}} 
			e^{-c|y-x|(t-s)^{-\epsilon/\alpha}}\int_{|u|\leq (t-s)^{\epsilon/\alpha}}  |u|^2 |\nu^x(du) - \nu^y(du)| \\
			& \lesssim (t-s)^{-\frac{d+2}{\alpha}} 
			e^{-c|y-x|(t-s)^{-\epsilon/\alpha}} (|x-y|^\eta \wedge 1) (t-s)^{\frac{(2-\alpha)\epsilon}{\alpha}}.
		\end{align*}
		Integrating in $y$, we derive $
		\sup_x\|I_1\|_{L_1} \lesssim (t-s)^{-1+\epsilon_1}$, $ 0<s<t\leq T$, 
		where $\epsilon_1:=\frac{\epsilon(d+\eta+2-\alpha)- (d+2-\alpha)}{\alpha}$. 
		Note that  we can chose $\epsilon$ close to 1; since $\eta\in (0,1]$, we have $\epsilon_1\in (0,1)$.

		Further, using \eqref{difx} with $k=0$ and  \eqref{compare},  we get
		$\sup_x	\|I_2\|_{L_1} \lesssim (t-s)^{-\epsilon}$, $ 0<s<t\leq T$. 
		Thus, 
		$$
		\sup_{x\in \rd} \|H_{s,t}(x,\cdot)\|_{L_1}\lesssim (t-s)^{-1+\delta_2}, \quad  0<s<t\leq T, 
		$$
		with $\delta_2:= \min(\epsilon_1, 1-\epsilon)$. 
		The same considerations are valid for $\tilde{H}$. 
		Thus, \eqref{HHeq}  holds true with $\delta_2$ as above. 
		
		Using the strategy justified in \cite{KK18}, \cite{KKS21}, one can show that i) the series \eqref{wt1} and \eqref{wt2} are convergent in $L_1$, ii)  $p$ and $\tilde{p}$ are transition probability densities of the processes $X$ and $\tilde{X}$, which are the unique solutions to the martingale problems for $(\LL, C_0^2)$ and $(\tilde{\LL}, C_0^2)$, respectively. 
		
		It remains to check \eqref{LLeq}.   Note that $
		(\widetilde{\mathcal L}_x-\mathcal L)_xp^0_{s,t}(x,y)=h^\kappa B_x p^0_{s,t}(x,y)$. 
		Using \eqref{difx}, we get 
		\begin{align*}
			\int_\rd|B_x p^0_{s,t}(x,y)|dy & \lesssim \int_\rd\int_{|u|\leq 1} \left(\int_{-1}^1 (1-|\theta|) |u' D^2_x p^0_{s,t}(y-x-\theta u) u| d\theta\right)  \frac{du}{|u|^{d+\beta}} dy \\
			& \quad + \int_\rd \left(\int_{|u|\geq 1} (p_{s,t}^0(y-x)+ p_{s,t}^0(x))\frac{du}{|u|^{d+\beta}} \right) dy\\
			& \lesssim (t-s)^{-\frac{d(1-\epsilon)}{\alpha}}
			\int_\rd (1\wedge (|u|(t-s)^{-1/\alpha})^2)|\frac{du}{|u|^{d+\beta}}  \lesssim (t-s)^{-\frac{d(1-\epsilon)+\beta}{\alpha}}. 
		\end{align*}
		Taking the $L_1$ norm in $y$ and using \eqref{compare}, we arrive at $ \sup_x	\|(\widetilde{\mathcal L}-\mathcal L)p^0_{s,t}(x,\cdot) \|_{L_1}\lesssim h^\kappa (t-s)^{-\frac{d(1-\epsilon)+\beta}{\alpha}}$. 
		Thus, \eqref{LLeq} is satisfied with $\delta_1 = 1-\frac{d(1-\epsilon)+\beta}{\alpha}$.   Since $\beta\in (0,\alpha)$ and $\epsilon$ can be taken close enough to 1, we can achieve that $\delta_1\in (0,1)$.    
	\end{proof}

	One can see the result of Theorem~\ref{T1} under a different angle. Namely, one can bound the weak error by  a certain distance between the respective characteristics $(b, 0, \nu)$ and $\widetilde b, 0,  \widetilde \nu)$. For the sake of simplicity, below we consider the case when both operators $\LL$ and $\widetilde \LL$ do not possess drifts.
	
	Let $d\geq 1$, $K: \rd \mapsto \rd\times \rd$; for $\gamma\in (0,1]$ define the (matrix) norms:
	$\|K\|_{\infty}:= \sup_{v\in \rd} \|K(v)\|_{op}$ and $[K]_{\CC^\gamma}:=  \sup_{0< |u|\leq 1, v\in \rd } \frac{\|K(v+u)- K(v)\|_{op}}{|u|^\gamma}$,  where $	\|A\|_{op}:= \sup_{|\ell|=1} |A\ell|$, $A \in \rd\times \rd$. 
	Define the $k$-the derivative of $K(v)$ and its operator norm in the following way:    For $k=0$ put $D^k K=K$,   and 
	$$
	D^k K(v)[\ell_1,\dot, \ell_k]:= \sum_{i_1, \dots i_k=1}^d (\ell_1)_{i_1}\dots (\ell_k)_{i_k}\frac{\partial^k K(v)}{\partial v_{i_1}\dots \partial v_{i_k}}, \quad k>0, 
	$$
	$$
	\|D^k K(v)\|:= \sup_{|\ell_1|=\dots |\ell_k|=1}\|D^k K(v)[\ell_1,\dot, \ell_k]\|_{op}, \quad \|D^k K\|_\infty:= \sup_{v\in \rd} \|D^k K(v)\|. 
	$$
	Define the (matrix) H\"older norms: 
	$	\|K\|_{\CC^\gamma}:= \sum_{k=0}^m \|D^k K\|_\infty + [D^m K]_{\CC^\lambda}$, where $\gamma = m +\lambda$, $m\in \mathbb{N} \cup \{0\}$, $\lambda \in (0,1)$. 
	If $K$ is a scalar function,  the above definition gives the ordinary H\"older space norm.

	Next, we define the   \emph{bounded $\gamma$-H\"older distance} between matrix-valued measures $m$ and $\widetilde m$: 
	\begin{equation}\label{BH}
		\dist_{BH,matr}^{(\gamma)}(m,\widetilde m):= \sup_{ \|K\|_{\CC^\gamma}\leq 1} \left| \int_\rd \tr \left(K (u)(m-\widetilde m)(du)\right)\right|, 
	\end{equation}
	where $K: \rd \mapsto \rd\times \rd$ is a matrix, and for scalar measures  $n$, $\widetilde n$
	\begin{equation}\label{BH2}
		\dist_{BH}^{(\gamma)}(n,\widetilde n):= \sup_{ \|\varphi\|_{\CC^\gamma}\leq 1} \left| \int_\rd \varphi(u)(n-\widetilde n)(du)\right|, 
	\end{equation}
	if $\varphi: \rd\mapsto \real$ is a scalar function.
	Note that we do not assume  the measures $m(du)$ and $\widetilde m(du)$  to have the same mass.  
	Let $S_r(u):= \frac{u}{r}$, $ \nu^x_r := (S_r)_\sharp \nu^x$, $\widetilde \nu^x_r := (S_r)_\sharp \widetilde \nu^x$,  and define
	$$
	n_r^x (dv):= \I_{|v|>1} \nu_r^x(dv),\qquad  \widetilde n_r^x (dv):=  \I_{|v|> 1}  \widetilde \nu_r^x(dv),$$  $$
	m_r^x (dv):= \I_{|v|\leq 1} vv'\nu_r^x(dv), \qquad \widetilde m_r^x (dv):=  \I_{|v|\leq 1} vv' \widetilde \nu_r^x(dv).
	$$
	Suppose that there exists $\beta \in (0,\alpha)$ such that
	\begin{equation}\label{Dgambet} 
		D_{\gamma, \beta}(\nu^x,
		\widetilde \nu^x):= \sup_{0<r\leq 1 } r^\beta \left(\dist_{BH,matr}^{(\gamma)}(m_r^x,\widetilde m_r^x) + \dist_{BH}^{(\gamma)}(n_r^x,\widetilde n_r^x)\right)<\infty. 
	\end{equation}
	
	\begin{theorem}\label{t3}
		Let $\LL$ and $\widetilde \LL$ be operators of type \eqref{L1}, with kernels $\nu^x$ and $\widetilde \nu^x$, respectively.  Let $\gamma\in (0,1)$, $\beta\in (0,\alpha)$ be such that \eqref{Dgambet}  holds true.  Suppose that the conditions \eqref{HHeq}  and 
		\begin{equation}\label{LLeq2}
			\sup_x\left\|
			(\widetilde \LL_x- \LL_x)p^0_{s,t}(x,\cdot)
			\right\|_{L_1}
			\lesssim  \sup_x D_{\gamma, \beta}(\nu^x,\widetilde  \nu^x)
			(t-s)^{-1+\delta_1}
		\end{equation}
		are satisfied. Then  for $f\in C_\infty(\rd)$ 
		\begin{equation}\label{te-eq1}
			\sup_x |\Ee^x f(X_t) - \Ee^x f(\widetilde X_t) | \lesssim t^{\delta}\|f\|_\infty\sup_x   D_{\gamma, \beta}(\nu^x,\widetilde  \nu^x), \quad 0<t\leq T, 
		\end{equation}
		where $\delta = \min(\delta_1,\delta_2)$. 
		In particular, 
		$$\sup_x \|\Pp_t(x,\cdot)  -\widetilde \Pp_t(x,\cdot) \|_{TV} \lesssim  t^\delta\sup_x   D_{\gamma, \beta}(\nu^x,\widetilde  \nu^x).
		$$
	\end{theorem}
	The proof follows  line by line the proof  of Theorem~\ref{T1}.  Consider  a simple case, when \eqref{LLeq2} holds true.   
	
	\begin{proposition}\label{Pr2}Suppose that $p^0$ satisfies \eqref{difx}. Then \eqref{LLeq2} holds true for any $\gamma \in (0,1)$. 
	\end{proposition}

	\begin{proof}[Proof of Proposition\eqref{Pr2}]
		We use the  $L_1-L_\infty$  duality relation:  $
		\|g\|_{L_1} = \sup_{\|\phi\|_\infty\leq 1} \left| \int_\rd \phi(y)g(y)dy\right|$ for any  $g\in L_1(\rd)$.
		Take $\phi \in C_\infty(\rd)$ and   $ g \equiv g_{s,t}(x,y):= (\LL_x - \widetilde \LL_x)p^0_{s,t}(x,y)$. 
		Using the Taylor expansion up to order 2,  Fubini's theorem,  splitting and the change of variables $u = vr$,  we get
		\begin{align*}
			\int_\rd \phi(y)& g_{s,t}(x,y)dy:= \frac12 \int_\rd \phi(y)\left[\int_{|u|\leq r} \left\{\int_{-1}^1 (1-|\theta|)u'D^2_x p_{s,t}^0 (x+\theta u,y)ud\theta \right\}(\nu^x  - \widetilde \nu^x)(du)\right]dy\\
			& \quad + \frac12 \int_\rd \phi(y)\left[\int_{|u|>r} \left(p_{s,t}^0(x+u,y) + p_{s,t}^0(x-u,y)-2p_{s,t}^0(x,y)\right)
			(\nu^x  - \widetilde \nu^x)(du)\right]dy\\
			& =\int_{\|v\|\leq 1}  \tr \left(K_{s,t}^{x,r,\phi}(v)  	(m_r^x - \widetilde m_r^x)(dv)\right) + \int_{\|v\|> 1}  F_{s,t}^{x,r,\phi}(v)  	(n_r^x - \widetilde n_r^x)(dv),
		\end{align*}
		where 
		\begin{equation}\label{K}
			K_{s,t}^{x,r,\phi}(v) := \frac{r^2}{2}\int_\rd \phi(y)\int_{-1}^1 (1-|\theta|) D^2_x p^0_{s,t}(x+\theta vr,y) d\theta dy
		\end{equation}
		and 
		\begin{equation}\label{F}
			F_{s,t}^{x,r,\phi}(v) :=
			\frac12 \int_\rd \phi(y) \left(p_{s,t}^0(x+rv,y) + p_{s,t}^0(x-rv,y)-2p_{s,t}(x,y)\right)dy. 
		\end{equation}
		%Here $\ell_v := v/|v|$ is the unit vector in the direction of $v$. 
		Let us estimate the  $\CC^\gamma$ -- norms of  $K_{s,t}^{x,r,\phi}(v)$ and $F_{s,t}^{x,r,\phi}(v)$. 
		Denote $	A_{s,t}^{(k)}:= \sup_x \|D_x^k p_{s,t}^{(0)}(x,\cdot)\|_{L_1}$,  $k\geq 0$. 
		Then 
		\begin{equation}\label{F1}
			\| D^k K_{s,t}^{x,r,\phi}\|_{C_\infty}\lesssim  r^{2+k} A_{s,t}^{(2+k)},  \quad k=0,1, 
		\end{equation}
		implying  
		\begin{equation}\label{F3}
			\begin{split}
				\| K_{s,t}^{x,r,\phi}\|_{\CC^\gamma}
				&\leq r^2 A_{s,t}^{(2)} + r^{2+\gamma} (A_{s,t}^{(2)})^{1-\gamma}(A_{s,t}^{(3)})^\gamma. 
			\end{split}
		\end{equation}
		where we used the inequality $\min(a,b) \leq a^{1-\gamma} b^{\gamma}$, $\gamma\in (0,1)$ for $a,b>0$. 
		
		Similarly, $\| D^k F_{s,t}^{x,r,\phi}\|_{C_\infty}\lesssim  r^k A_{s,t}^{(k)}$,  $ k=0,1$, 
		implying 
		\begin{equation}\label{F4}
			\begin{split}
				\| F_{s,t}^{x,r,\phi}\|_{\CC^\gamma}\leq A_{s,t}^{(0)} + r^{\gamma} (A_{s,t}^{(0)})^{1-\gamma}(A_{s,t}^{(1)})^\gamma. 
			\end{split}
		\end{equation}
		Taking now  $r= (t-s)^{\epsilon/\alpha}$, we arrive at 
		\begin{equation}\label{F5}
			\begin{split}
				\left|	\int_\rd  \phi(y)g_{s,t}(x,y)dy\right|&\lesssim  (t-s)^{-\frac{(d+2+\gamma)(1-\epsilon)}{\alpha}}\left( \dist_{BH}^{(\gamma)} (m_r^x, \widetilde m_r^x) + \dist_{BH}^{(\gamma)} (m_r^x, \widetilde m_r^x)\right)\\
				& \lesssim  (t-s)^{-\frac{(d+2+\gamma)(1-\epsilon)}{\alpha}- \frac{\epsilon \beta}{\alpha}}\sup_{r\in [0,1]}r^{\beta}\left( \dist_{BH}^{(\gamma)} (m_r^x, \widetilde m_r^x) + \dist_{BH}^{(\gamma)} (m_r^x, \widetilde m_r^x)\right)\\
				& = (t-s)^{-1+ \delta_1}\sup_x  D_{\gamma, \beta}(\nu^x, \widetilde \nu^x), 
			\end{split}
		\end{equation}
		where 
		$\delta_1 := 1- \frac{(d+2+\gamma)(1-\epsilon)+\beta \epsilon}{\alpha}$. Since $\beta\in (0,\alpha)$, we have $\delta_1>0$ provided that we choose $\epsilon> \frac{d+2+\gamma-\alpha}{d+2+\gamma-\beta}$. 
	\end{proof}
	The next example illustrates that  condition \eqref{LLeq2} can be satisfied also for L\'evy(type) measures  which are not mutually absolutely continuous. 
	\begin{example}\label{exa2}
		Consider the one-dimensional case.   For  $0<\beta<\alpha$, define the L\'evy measure
		\begin{equation}\label{nub}
			\nu^{(\beta)}(du):= \sum_{n=1}^\infty 2^{n\beta}\left(\delta_{2^{-n}}+ \delta_{-2^{-n}}\right)(du), 
		\end{equation}
		the transform $T_h(u):= (1+h)u$, $h\in (0,1)$, 
		and the push-forward L\'evy measure
		\begin{equation}\label{nub} 
			\nu^{(\beta,h)}(du):= (T_h)_\sharp \nu^{(\beta)}(du)= \sum_{n=1}^\infty 2^{n\beta}\left(\delta_{2^{-n}(1+h)}+ \delta_{-2^{-n}(1+h)}\right)(du). 
		\end{equation}
		Consider the following operators: $
		\Aa f(x) = \LL f(x) + B^{(\beta)}f$,  $\widetilde \Aa f(x) = \LL f(x) + B^{(\beta,h)}f$, 
		where $\LL$ is defined in \eqref{L1}, $B^{(\beta)}$ and $B^{(\beta,h)}$  are symmetric L\'evy operators with zero drift and diffusion components, and L\'evy measures $	\nu^{(\beta,h)}$ and 	$\nu^{(\beta,h)}$, respectively. 
		Note that for $h\in (0,1)$ the measures   $	\nu^{(\beta)}$ and 	$\nu^{(\beta,h)}$ are mutually singular. Indeed, for  $2^{-n}= (1+h)2^{-m}$, we should have $1+h = 2^{n-m}$, which is impossible for $h\in (0,1)$.

		Relying on Proposition~\ref{Pr1}, one can construct the 
		parametrix series for each of the operators, starting with the same zero-order approximation $p^0$, and show that  \eqref{HHeq} is satisfied.  Further, using \eqref{difx}, one can apply Proposition~\ref{Pr2} and get \eqref{LLeq2}.  The latter is satisfied, if  for $\psi\in \CC^{\gamma}$ with $\|\psi\|_{\CC^{\gamma}}\leq 1$,  $\gamma \in (0,1)$, we show that for 
		$\mu^{(\beta)}(dv) = v^2 \nu^{(\beta)}(dv)$, $ \widetilde \mu^{(\beta,h)}(dv) = v^2 \widetilde\nu^{(\beta,h)}(dv)$, 
		one has
		\begin{equation}\label{exa2-eq20}
			\left|\int_\real \psi(u)(\mu^{(\beta)} -\widetilde \mu^{(\beta,h)})(du)\right|\lesssim  h^{\gamma}. 
		\end{equation} 
		Indeed, from 
		$$
		\int_\real \psi(v)\widetilde \mu^{(\beta,h)}(dv) =\int_\real T_h^2(u) \psi(T_h(u))\nu^{(\beta)} (du), \qquad 
		\int_\real \psi(v) \mu^{(\beta)}(dv) =\int_\real u^2 \psi(u)\nu^{(\beta)} (du), 
		$$
		we get, using $\|\psi\|_{\CC^{\gamma}}\leq 1$, 
		\begin{align*}
			\left|\int_\real \psi(u)(\mu^{(\beta)} -\widetilde \mu^{(\beta,h)})(du)\right| & \leq \int_\real u^2 \left| \psi(u)- \psi(T_h(u))\right|\nu^{(\beta)} (du)+ \int_\real \left|(u^2 -T^2_h(u))\psi(T_h(u))\right|\nu^{(\beta)}(du)\\
			%	&\leq \int_\real u^2 \left| \psi(u)- \psi(T_h(u))\right|\nu^{(\beta)} (du)+ \int_\real \left|(u^2 -T^2_h(u))\right|\nu^{(\beta)}(du)\\\
			& \lesssim h^{\gamma }\int_\real |u|^{2+\gamma}\nu^{(\beta)}(du) + (2h +h^{2})\int_\real |u|^2 \nu^{(\beta)}(du) \lesssim h^{\gamma}, 
		\end{align*}
		justifying \eqref{exa2-eq20}. 
	\end{example}
	
	\begin{remark}
		In this section we made a relatively strong assumption on the existence of the transition probability densities for both processes $X$ and $\tilde X$, their parametrix representation and, moreover, that the zero order approximations in both series coincide. This approach has an advantage when the zero order approximation  can be chosen in such a way that it has ``good properties'' (i.e., its gradient estimates are suitable for the proof of convergence). On the other hand, it is obviously a limitation. In the next section we set aside this assumption  and prove the result, in a sense similar to that of Theorem~\ref{t3}, under the assumption of  existence of the transition probability density for one of the processes  its gradient estimates. 
	\end{remark}

	\subsection{Weak error II:  small-jump Gaussian approximation}
	
	In this section for the sake of transparency we assume that assumptions \eqref{compare} and \eqref{Lip} on the kernel $\nu^x(\cdot)$ are satisfied. 
	We do not assume that the kernel $\nu^x(\cdot)$ is symmetric. 
	
	Consider $\LL$ of type \eqref{LL}; 
	for $T>0$, $\epsilon,\delta\in [0,1)$,  and the cutting level $	r_{\epsilon, \delta}(t):= \delta t^{\epsilon/\alpha}$, 
	define 
	\begin{equation}\label{tL}
		\begin{split}
			\LL^{\epsilon,\delta} f(x)
			& :=b(x)\cdot D f(x)+\frac12  
			\int_{|u|\le r_{\epsilon, \delta}(T-t)}
			u'D^2 f(x)u \nu^x(du)\\
			&\quad +\int_{|u|> r_{\epsilon, \delta}(T-t)}
			\bigl(f(x+u)-f(x)-D f(x)\cdot u \I_{|u|\leq 1}\bigr)\,\nu^x(du),
			\quad 0<t<T. 
		\end{split}
	\end{equation}
	%where 
	%$$
	%b(t,x) = b(x) - %\int_{r_{\epsilon,\delta}(T-t)<|u|\leq 1} u \nu(x,du).
	%$$ 
	%Let $C_\infty(\rd)$ be the space of continuous functions vanishing at infinity with the $\sup$-norm $\|\cdot\|_\infty$ on it.

	Assume that $\delta T^{\epsilon/\alpha} \leq 1$, otherwise we need  to change appropriately the drift function $b(x)$. 	Further, let $\widetilde \LL$ be defined in \eqref{tL}.	
	Suppose that $\lambda, \eta>0$ satisfy the restrictions 
	\begin{equation}\label{As3eq1}
		(2+\lambda - \alpha)(1-\epsilon)<\eta < 2+ \lambda  
	\end{equation} 
	and that  the conditions below are satisfied: 
	
	\begin{itemize} 
		
		\item[\textbf{(V1)}] For any $f\in \CC^\eta (\rd)$ for $\eta$ satisfying \eqref{As3eq1}, the function $v(s,x):= \Ee^x f(X_{T-s})$ is the unique   bounded solution  to the Cauchy problem 
		\begin{equation}\label{Cauch}
			(\partial_s  + \LL_x) v(s,x)=0,
			\qquad
			v(T,x)=f(x); 
		\end{equation}
		%	where $f\in \CC^\eta (\rd)$ for some $\eta \in \real$; 
		
		\item[\textbf{(V2)}] For $\lambda\in (0,1]$ satisfying \eqref{As3eq1}, 
		$$
		\| v(s,\cdot)\|_{C^{2+\lambda}}\lesssim (T-s)^{-(2+\lambda - \eta)/\alpha}\|f\|_{\CC^\eta}, \quad 0<s<T. 
		$$
		
		%		\item[\textbf{(V3)}] \red $\lim_{s\uparrow T} v(s,x)= f(x)$, locally uniformly in $x\in \rd$. \normal
	\end{itemize}
	
	\begin{theorem}\label{t2}
		Let $f\in \CC^\eta(\rd)$, where $\eta$ satisfies \eqref{As3eq1}.  Under Assumptions \textbf{(V1)} and \textbf{(V2)}, 
		\begin{equation}\label{t2eq1}
			\sup_x |\Ee^x f(X_T) - \Ee^x f(X_T^{\epsilon,\delta}) |\leq C_T \delta^{2+\lambda -\alpha} \|f\|_{\CC^\eta},
		\end{equation} 
		where $C_T \lesssim  \frac{\alpha T^{\frac{\eta - (2+\lambda - \alpha)(1-\epsilon)}{\alpha}}}{\eta - (2+\lambda - \alpha)(1-\epsilon)}$. 
	\end{theorem}

	%Before we proceed to the proof, let us discuss the distinction between the result in the case of constant-level  (i.e.m  when $\epsilon=0$) and  dynamic cuttings. 
	
	\begin{remark}\label{rem2}
		Note that both ways -- dynamic and constant  cutting -- give the same order of decay  in $\delta$ as $\delta \to 0$ in the  right-hand side in \eqref{t2eq1}. However, making the dynamic cutting allows to take less regular functions $f$. Namely, making the constant level cutting, we   derive \eqref{t2eq1} with $f\in \CC^\eta$ where $\eta > 2+\lambda - \alpha$, whereas the dynamic cutting allows to take $\eta$ smaller, according to \eqref{As3eq1}. Note that this restriction becomes stronger when $\alpha$ is close to 0. 
	\end{remark}
	
	\begin{proof}[Proof of Theorem~\ref{t2}]
		Using the assumption  $X_0=X_0^{\epsilon,\delta} =x$,  we can apply the Ito formula on $[0,T)$ to $v(s,X_s)$ and $v(s,X_s^{\epsilon,\delta})$, respectively: For $t\in (0,T)$, 
		\begin{equation}\label{t2-eq1}
			\begin{split}
				v(t,X_t)- v(0,x) & = \int_0^t (\partial_s + \LL_x)v(s,X_s)ds + M_t, 
				\\
				v(t,X_t^{\epsilon,\delta})- v(0,x)& = \int_0^t (\partial_s + \LL^{\epsilon,\delta}_x)v(s,X_s^{\epsilon,\delta})ds + \tilde{M}_t
			\end{split}
		\end{equation}
		where $(M_t)_{t\in [0,T)}$, $(\tilde{M}_t)_{t\in [0,T)}$,  are  local martingales on $[0,T)$.  By \textbf{(V1)}, \textbf{(V2)} and the uniform boundedness of coefficients, assumed through  the paper,  both local martingales are bounded on $[0,t]$, $t<T$, and hence are true martingales on this interval.  By assumption \textbf{(V1)}, 
		$\partial_s v = - \LL_x v$, hence, 
		$$
		(\partial_s + \LL^{\epsilon,\delta}_z) v(s,z)= (\LL^{\epsilon,\delta}_z- \LL_z) v(s,z). 
		$$
		
		Subtracting and taking the expectation gives 
		\begin{equation}\label{t2-eq2}
			\Ee^x v(t,X_t^{\epsilon,\delta})-\Ee^x v(t,X_t)
			=
			\Ee^x\int_0^t
			(\LL^{\epsilon,\delta}_x- \LL_x)
			v(s, X_s^{\epsilon,\delta})\,ds.
		\end{equation}
		Since the large-jump parts of $\LL$ and $\LL^{\epsilon,\delta}$ coincide,
		\begin{align*}
			(\LL_z-\LL^{\epsilon,\delta}_z)v(s,z)
			=&
			\int_{|u|\le r_{\epsilon,\delta}(T-s)}
			\biggl[
			v(s,z+u)-v(s,z)-D_z v(s,z)\cdot u
			%	\\
			%	&\hspace{35mm}
			-\frac12u' D^2_z v(s,z) u 
			\biggr]\nu(z,du).
		\end{align*}
		Using \textbf{(V2)}, we can estimate the inner expression as 
		\begin{align*}
			\big|	v(s,z+u)-v(s,z)&-D_zv(s,z)\cdot u
			-\frac12u' D^2_z v(s,z) u \big|\\
			&=\left| \int_0^1 (1-\theta)u' [D^2_z v(s,z+\theta u)- D^2_z v(s,z)]u \, d\theta\right|\\
			& \lesssim \|v(s,\cdot)\|_{\CC^{2+\lambda}} |u|^{2+\lambda}. 
		\end{align*}
		Substitution the above estimate in \eqref{t2-eq2} and taking expectation  gives
		\[
		\begin{aligned}
			\left|\Ee^x v(t,X_t^{\epsilon,\delta})-\Ee^x v(t,X_t)
			\right|
			&\lesssim \|f\|_{\CC^\eta}
			\int_0^t
			(T-s)^{-(2+\lambda - \eta)/\alpha} \left(\sup_z\int_{|u|\le r_{\epsilon,\delta}(T-s)}|u|^{2+\lambda}
			\nu(z,du)
			\right)ds\\
			& \lesssim  \delta^{2+\lambda - \alpha} \frac{\alpha T^{\frac{\eta  - (2+\lambda - \alpha)(1-\epsilon)}{\alpha}}}{\eta- (2+\lambda - \alpha)(1-\epsilon)}  \|f\|_{\CC^\eta}, 
		\end{aligned}
		\]
		provided that \eqref{As3eq1} is satisfied.  
		Note that  $v(t,x)= P_{T-t}f(x):=\Ee^x f(X_{T-t})$. Hence,  by the Markov property we have 
		$$
		\Ee^x v(t,X_t) = P_t P_{T-t}f(x) = P_T f(x) = \Ee^x f(X_T). 
		$$
		It remains to show that 
		\begin{equation}\label{t2-eq3}
			\lim_{t\uparrow T} \Ee^x P_{T-t} f(X_t^{\epsilon, \delta}) = \Ee^x f(X_T^{\epsilon,\delta}). 
		\end{equation}
		By our standing assumption, $X_t^{\epsilon,\delta}$ is c\'adl\'ag, hence, $X_t^{\epsilon,\delta} \rightarrow X_{T-}^{\epsilon,\delta}$ as $t\uparrow T$. Note that by \textbf{(N1)}  for every $a>0$ the intensity of jumps of $X^{\epsilon,\delta}$ of  size greater than $a$ is uniformly bounded. Therefore, the probability of such a jump occurring in a time integral of length $h$ is less than $C_a h$, and tends to 0 as $h\to 0$. Thus, $X_t^{\epsilon,\delta}$ has no jump any  fixed deterministic  time moment, in particular, $X_{T-}^{\epsilon,\delta} = X_T^{\epsilon,\delta}$  $\Pp$--a.s.
		Note that by the assumption that the coefficients of \eqref{tL} are bounded, we get (cf. \cite[Th.6.4]{KSch23}, see also \cite{Sch98})
		\begin{equation}\label{t2-eq5}
			\lim_{h\downarrow 0}\sup_{x\in \rd}\Pp^x (|X_h-x|\geq \rho)=0 \quad \text{for any $\rho>0$.}
		\end{equation}
		Hence,  $P_{T-t}f(x)\to f(x)$  locally uniformly as $t\uparrow T$ for any $f\in \CC^\eta$. For $R>0$ large enough split 
		\begin{align*}
			\Ee^x |P_{T-t} f(X_t^{\epsilon,\delta}) - f(X_T^{\epsilon,\delta})| & \leq \sup_{|z|\leq R} |P_{T-t} f(z)-f(z)| \quad + 2\|f\|_\infty \Pp^x (|X_t^{\epsilon,\delta}|>R)+ \Ee^x |f(X_t^{\epsilon,\delta})- f(X_T^{\epsilon,\delta})| \\
			&= : I_1+I_2+I_3. 
		\end{align*}
		The first term tend to 0 as $t\uparrow T$ by  locally uniform convergence, the third one -- by the dominated convergence theorem. For the second term, since by the c\'adl\'ag property  the family of random variables $\{X_t^{\epsilon, \delta}:\, t\in [T/2,T]\}$ is tight, we have $
		\lim_{R\to \infty} \sup_{t\in [T/2,T]}\Pp^x (|X_t^{\epsilon,\delta}|>R)=0$. 
		Hence, the second term can be made  arbitrary small (uniformly in $t\in [T/2,T]$) by choosing $R$ large enough. Thus, we get  \eqref{t2-eq3}, which completes the proof. 
	\end{proof}

	\subsection{Examples}
	We begin with the space-homogeneous case. 
	\begin{example}\label{2ex1}
		Take  $\LL = - (-\bigtriangleup)^{\alpha/2}$; then $\nu(du) = c_{d,\alpha} |u|^{-d-\alpha}du$ for some positive constant $c_{d,\alpha}$.   The respective ``gaussianized'' operator $\LL^{\epsilon,\delta}$ is then given by 
		\begin{equation}\label{2exeq1}
			\LL^{\epsilon,\delta}f(x) = \int_{|u|> r_{\epsilon,\delta}(T-t)}(f(x+u)-f(x))
			\frac{c_{d,\alpha}}{|u|^{d+\alpha}}du +  \frac{Q^{\epsilon, \delta}(T-t)}{2}\bigtriangleup f (x). 
		\end{equation}
		where $Q^{\epsilon,\delta}(t):= \frac{c_{d,\alpha}\Vol(S^{d-1})}{d(2-\alpha)}\delta^{2-\alpha}t^\frac{(2-\alpha)\epsilon}{\alpha}$, where $\Vol(S^{d-1})$ is the volume of the unit  sphere in $\rd$.   Let  $v(s,x)= \Ee^x f(X_{T-s})$, $f\in \CC^\eta$. Then  condition \textbf{(V2)} follows from the ``caloric smoothing''  estimate proved in \cite[Thm.4.1]{KS22},  taking  the Besov spaces $B^{2+\lambda}_{\infty, \infty}$ and $B^{\eta}_{\infty, \infty}$, respectively (i.e., taking the Besov spaces of order  $d+s = 2+\lambda$ and  $s= \eta$, respectively, and recalling that for $s>0$ the Besov space $B^{s}_{\infty,\infty}$ coincides with the H\"older space $\CC^s$). Then $v(s,\cdot)\in C^2_\infty$, and thus \textbf{(V1)} is satisfied. 
	\end{example}
	\begin{example}\label{2ex2}
		The above  example can be easily generalized to the case $\LL  f= - \phi(- \bigtriangleup)$, where $\phi$ is a Bernstein function, $\phi(0)=0$, without the linear term,   satisfying $\phi(t) \geq C t^{\alpha/2}$ for $t\geq 1$, and 
		$\LL^{\epsilon,\delta}$ of the form \eqref{2exeq1} with $Q^{\epsilon,\delta}(t) =\frac{1}{d} \int_{|u|\leq r_{\epsilon,\delta}(t)} |u|^2 \nu(du)$, where $\nu(\cdot)$ is the L\'evy measure corresponding to $\LL$. Note that under the restriction $\epsilon<\frac{\alpha}{2}$ the process, corresponding to the L\'evy measure $\nu(du) \I_{|u|\geq r_{\epsilon,\delta}(t)}$, is an (inhomogeneous)  compound Poisson with integral  intensity 
		\begin{align*}
			\Lambda(T) & := \int_0^T\int_{|u|> r_{\epsilon,\delta}(t)} \nu(du) dt \lesssim \frac{\delta^{-2} T^{1-2\epsilon/\alpha}}{1-2\epsilon/\alpha}<\infty. 
		\end{align*}
		In the above calculation  we used that 
		$\nu(|u|>r)\leq \nu(|u|>1) + r^{-2}\int_{r<|u|\leq 1} |u|^2\nu(du)\lesssim r^{-2}$. 
		Thus,  using the cutting level $r_{\epsilon, \delta}$, in order to have a simplified model of type  ``Compound Poisson process + Diffusion'',  it is sufficient  to assume $\epsilon<\frac{\alpha}{2}$. Condition \textbf{(V2)} is provided by \cite[Cor.5.2]{KS22}, since $(\phi^{-1}(1/t))^d\lesssim t^{-d/\alpha}$, $t\in (0,T]$. Moreover, the semigroup representation of $v(s,x)$ shows that   $v(s,x)$ is the classical solution to  \eqref{Cauch}; thus \textbf{(V1)} holds true. 
	\end{example}
	In the next example we consider a L\'evy-type process and  provide the conditions and references, under which   \textbf{(V1)},  \textbf{(V2)}  are satisfied. 
	
	\begin{example}
		We quote the result from \cite{CHZ20} but under our restriction \eqref{As3eq1} on $\eta$  (the initial result, stated in \cite{CHZ20},  holds for more general set of parameters).  Let
		\begin{equation}
			\LL f(x) = \sum_{k=1}^d p.v. \int_\real (f(x+\sigma^{\cdot k}(x) z )- f(x))\frac{dz}{|z|^{1+\alpha}} + b(x)\cdot D f(x), 
		\end{equation}
		where $\sigma(x):  \rd\mapsto \rd\otimes \rd$ is a matrix-valued measurable function, satisfying the uniform ellipticity condition: There exists $c_0>0$, s.t.
		$c_0^{-1} |\xi| \leq |\sigma(x)\xi| \leq c_0 |\xi|$, $ \forall t\geq 0$, $x,\xi \in \rd$.

		Assume also that  $\| \sigma(\cdot)\|_{\CC^1}\lesssim 1$, and one of the conditions below is satisfied: 
		
		$$
		a) \quad b=0,\quad \beta =1, \qquad \text{or} \qquad b) \quad \alpha\in (1/2,2), \quad  \|b(\cdot)\|_{\CC^\beta}\lesssim 1, 
		$$
		where $\alpha+\beta>1$. Then   the martingale problem for $(\LL, C_0^2)$ is well-posed (see \cite{CHZ20}).
		Moreover, (cf. \cite[(1.4)]{CHZ20}),  for any $\lambda, \eta$ such that $0<\lambda< \min\{1,\alpha+ \min(\alpha,\beta)-2\}$, $(2+\lambda-\alpha)(1-\epsilon)<\eta<2+\lambda$, the estimated stated in \textbf{(V2)} holds true. By the same argument as in the previous examples, then   \textbf{(V1)}  is satisfied. 
		
	\end{example}

	\bibliographystyle{abbrv}
	
	\bibliography{prmx-diff}

@article {KSch23,
    AUTHOR = {K\"{u}hn, Franziska and Schilling, Ren\'{e} L.},
     TITLE = {Maximal inequalities and some applications},
   JOURNAL = {Probab. Surv.},
  FJOURNAL = {Probability Surveys},
    VOLUME = {20},
      YEAR = {2023},
     PAGES = {382--485},
      ISSN = {1549-5787},
   MRCLASS = {60E15 (42A61 42B25 60G15 60G44 60G51 60G53 60J25)},
  MRNUMBER = {4568704},
       DOI = {10.1214/23-ps17},
       URL = {https://doi.org/10.1214/23-ps17},
}

@article {Sch98,
    AUTHOR = {Schilling, Ren\'{e} L.},
     TITLE = {Growth and {H}\"{o}lder conditions for the sample paths of
              {F}eller processes},
   JOURNAL = {Probab. Theory Related Fields},
  FJOURNAL = {Probability Theory and Related Fields},
    VOLUME = {112},
      YEAR = {1998},
    NUMBER = {4},
     PAGES = {565--611},
      ISSN = {0178-8051,1432-2064},
   MRCLASS = {60J75 (35S99 47D07 60F15 60G17)},
  MRNUMBER = {1664705},
MRREVIEWER = {Niels\ Jacob},
       DOI = {10.1007/s004400050201},
       URL = {https://doi.org/10.1007/s004400050201},
}

@article {GHK18,
    AUTHOR = {Gairing, Jan and H\"{o}gele, Michael and Kosenkova, Tetiana},
     TITLE = {Transportation distances and noise sensitivity of
              multiplicative {L}\'{e}vy {SDE} with applications},
   JOURNAL = {Stochastic Process. Appl.},
  FJOURNAL = {Stochastic Processes and their Applications},
    VOLUME = {128},
      YEAR = {2018},
    NUMBER = {7},
     PAGES = {2153--2178},
      ISSN = {0304-4149,1879-209X},
   MRCLASS = {60G51 (60G52 60H10 60J75 62M10 62P12)},
  MRNUMBER = {3804790},
       DOI = {10.1016/j.spa.2017.09.003},
       URL = {https://doi.org/10.1016/j.spa.2017.09.003},
}

@article {GHKK15,
    AUTHOR = {Gairing, Jan and H\"{o}gele, Michael and Kosenkova, Tetiana
              and Kulik, Alexei},
     TITLE = {Coupling distances between {L}\'{e}vy measures and
              applications to noise sensitivity of {SDE}},
   JOURNAL = {Stoch. Dyn.},
  FJOURNAL = {Stochastics and Dynamics},
    VOLUME = {15},
      YEAR = {2015},
    NUMBER = {2},
     PAGES = {1550009, 25},
      ISSN = {0219-4937,1793-6799},
   MRCLASS = {60J60 (60F17 60G51 60H10 60J75)},
  MRNUMBER = {3332269},
MRREVIEWER = {Yulin\ Song},
       DOI = {10.1142/S0219493715500094},
       URL = {https://doi.org/10.1142/S0219493715500094},
}

@article {KK14,
    AUTHOR = {Kosenkova, T. and Kulik, A.},
     TITLE = {Transportation distance between the {L}\'{e}vy measures and
              stochastic equations for {L}\'{e}vy-type processes},
   JOURNAL = {Mod. Stoch. Theory Appl.},
  FJOURNAL = {Modern Stochastics. Theory and Applications},
    VOLUME = {1},
      YEAR = {2014},
    NUMBER = {1},
     PAGES = {49--64},
      ISSN = {2351-6046,2351-6054},
   MRCLASS = {60J75 (60G51 60H10)},
  MRNUMBER = {3314793},
MRREVIEWER = {Hannah\ Geiss},
       DOI = {10.15559/vmsta-2014.1.1.7},
       URL = {https://doi.org/10.15559/vmsta-2014.1.1.7},
}

@misc{Hu16,
      title={Density stability for some L\'evy-driven Stochastic Differential Equations}, 
      author={L Huang},
      year={2016},
      eprint={1603.05088},
      archivePrefix={arXiv},
      primaryClass={math.PR},
      url={https://arxiv.org/abs/1603.05088}, 
}

@article {KM11,
    AUTHOR = {Konakov, Valentin and Menozzi, St\'{e}phane},
     TITLE = {Weak error for stable driven stochastic differential
              equations: expansion of the densities},
   JOURNAL = {J. Theoret. Probab.},
  FJOURNAL = {Journal of Theoretical Probability},
    VOLUME = {24},
      YEAR = {2011},
    NUMBER = {2},
     PAGES = {454--478},
      ISSN = {0894-9840,1572-9230},
   MRCLASS = {60H10 (60G52 65C30)},
  MRNUMBER = {2795049},
MRREVIEWER = {Robert\ Wooster},
       DOI = {10.1007/s10959-010-0291-x},
       URL = {https://doi.org/10.1007/s10959-010-0291-x},
}

@article {KKM17,
    AUTHOR = {Konakov, Valentin and Kozhina, Anna and Menozzi, St\'{e}phane},
     TITLE = {Stability of densities for perturbed diffusions and {M}arkov
              chains},
   JOURNAL = {ESAIM Probab. Stat.},
  FJOURNAL = {ESAIM. Probability and Statistics},
    VOLUME = {21},
      YEAR = {2017},
     PAGES = {88--112},
      ISSN = {1292-8100,1262-3318},
   MRCLASS = {60H10 (60J60 65C30)},
  MRNUMBER = {3630604},
MRREVIEWER = {Yuguo\ Lin},
       DOI = {10.1051/ps/2016028},
       URL = {https://doi.org/10.1051/ps/2016028},
}

@article {BK21,
    AUTHOR = {Bitter, Ilya and Konakov, Valentin},
     TITLE = {{$L_1$} and {$L_\infty$} stability of transition densities of
              perturbed diffusions},
   JOURNAL = {Random Oper. Stoch. Equ.},
  FJOURNAL = {Random Operators and Stochastic Equations},
    VOLUME = {29},
      YEAR = {2021},
    NUMBER = {4},
     PAGES = {287--308},
      ISSN = {0926-6364,1569-397X},
   MRCLASS = {60H20 (60G22 60H15)},
  MRNUMBER = {4343186},
       DOI = {10.1515/rose-2021-2067},
       URL = {https://doi.org/10.1515/rose-2021-2067},
}

@book {Fr64,
    AUTHOR = {Friedman, Avner},
     TITLE = {Partial differential equations of parabolic type},
 PUBLISHER = {Prentice-Hall, Inc., Englewood Cliffs, NJ},
      YEAR = {1964},
     PAGES = {xiv+347},
   MRCLASS = {35.00 (35.62)},
  MRNUMBER = {181836},
MRREVIEWER = {B.\ Frank\ Jones, Jr.},
}

@article {KKS21,
    AUTHOR = {Knopova, Victoria and Kulik, Alexei and Schilling, Ren\'{e}
              L.},
     TITLE = {Construction and heat kernel estimates of general stable-like
              {M}arkov processes},
   JOURNAL = {Dissertationes Math.},
  FJOURNAL = {Dissertationes Mathematicae},
    VOLUME = {569},
      YEAR = {2021},
     PAGES = {86},
      ISSN = {0012-3862,1730-6310},
   MRCLASS = {60J35 (35A08 35A17 35S05 60G52 60J25)},
  MRNUMBER = {4361582},
MRREVIEWER = {A.\ Entsar\ Abdel-Rehim},
       DOI = {10.4064/dm824-8-2021},
       URL = {https://doi.org/10.4064/dm824-8-2021},
}

@article {KK18,
    AUTHOR = {Knopova, Victoria and Kulik, Alexei},
     TITLE = {Parametrix construction of the transition probability density
              of the solution to an {SDE} driven by {$\alpha$}-stable noise},
   JOURNAL = {Ann. Inst. Henri Poincar\'{e} Probab. Stat.},
  FJOURNAL = {Annales de l'Institut Henri Poincar\'{e} Probabilit\'{e}s et
              Statistiques},
    VOLUME = {54},
      YEAR = {2018},
    NUMBER = {1},
     PAGES = {100--140},
      ISSN = {0246-0203,1778-7017},
   MRCLASS = {60J35 (35R11 35S05 35S10 47G30 60J75)},
  MRNUMBER = {3765882},
MRREVIEWER = {Arturo\ Kohatsu-Higa},
       DOI = {10.1214/16-AIHP796},
       URL = {https://doi.org/10.1214/16-AIHP796},
}

@article {KK13,
    AUTHOR = {Knopova, Victoria and Kulik, Alexei},
     TITLE = {Intrinsic small time estimates for distribution densities of
              {L}\'{e}vy processes},
   JOURNAL = {Random Oper. Stoch. Equ.},
  FJOURNAL = {Random Operators and Stochastic Equations},
    VOLUME = {21},
      YEAR = {2013},
    NUMBER = {4},
     PAGES = {321--344},
      ISSN = {0926-6364,1569-397X},
   MRCLASS = {60G51 (41A60 60J75)},
  MRNUMBER = {3139314},
MRREVIEWER = {Ying\ Hui\ Dong},
       DOI = {10.1515/rose-2013-0015},
       URL = {https://doi.org/10.1515/rose-2013-0015},
}

@article {AR01,
    AUTHOR = {Asmussen, Soeren and Rosi\'{n}ski, Jan},
     TITLE = {Approximations of small jumps of {L}\'{e}vy processes with a
              view towards simulation},
   JOURNAL = {J. Appl. Probab.},
  FJOURNAL = {Journal of Applied Probability},
    VOLUME = {38},
      YEAR = {2001},
    NUMBER = {2},
     PAGES = {482--493},
      ISSN = {0021-9002,1475-6072},
   MRCLASS = {60G51 (60E07 60F05 60F17 60G52)},
  MRNUMBER = {1834755},
MRREVIEWER = {Ornella\ Arcudi},
       DOI = {10.1017/s0021900200019987},
       URL = {https://doi.org/10.1017/s0021900200019987},
}

@article {BM25,
    AUTHOR = {Bossy, Mireille and Maurer, Paul},
     TITLE = {On the {$\varepsilon$}-{E}uler-{M}aruyama scheme for time
              inhomogeneous jump-driven {SDE}s},
   JOURNAL = {Stochastic Process. Appl.},
  FJOURNAL = {Stochastic Processes and their Applications},
    VOLUME = {190},
      YEAR = {2025},
     PAGES = {Paper No. 104747, 34},
      ISSN = {0304-4149,1879-209X},
   MRCLASS = {60H35 (60H10 60K50)},
  MRNUMBER = {4943239},
MRREVIEWER = {Annika\ Lang},
       DOI = {10.1016/j.spa.2025.104747},
       URL = {https://doi.org/10.1016/j.spa.2025.104747},
}

@article {DMT21,
    AUTHOR = {Deligiannidis, G. and Maurer, S. and Tretyakov, M. V.},
     TITLE = {Random walk algorithm for the {D}irichlet problem for
              parabolic integro-differential equation},
   JOURNAL = {BIT},
  FJOURNAL = {BIT. Numerical Mathematics},
    VOLUME = {61},
      YEAR = {2021},
    NUMBER = {4},
     PAGES = {1223--1269},
      ISSN = {0006-3835,1572-9125},
   MRCLASS = {65M75 (35R09 60H10 60H35 60J74)},
  MRNUMBER = {4332136},
       DOI = {10.1007/s10543-021-00863-2},
       URL = {https://doi.org/10.1007/s10543-021-00863-2},
}

@article {CHZ20,
    AUTHOR = {Chen, Zhen-Qing and Hao, Zimo and Zhang, Xicheng},
     TITLE = {H\"{o}lder regularity and gradient estimates for {SDE}s driven
              by cylindrical {$\alpha$}-stable processes},
   JOURNAL = {Electron. J. Probab.},
  FJOURNAL = {Electronic Journal of Probability},
    VOLUME = {25},
      YEAR = {2020},
     PAGES = {Paper No. 137, 23},
      ISSN = {1083-6489},
   MRCLASS = {60H10 (60G52)},
  MRNUMBER = {4179301},
MRREVIEWER = {Josep\ Vives},
       DOI = {10.1214/20-ejp542},
       URL = {https://doi.org/10.1214/20-ejp542},
}

@article {KS22,
    AUTHOR = {Knopova, Victoriya and Schilling, Ren\'{e} L.},
     TITLE = {Bochner's subordination and fractional caloric smoothing in
              {B}esov and {T}riebel-{L}izorkin spaces},
   JOURNAL = {Math. Nachr.},
  FJOURNAL = {Mathematische Nachrichten},
    VOLUME = {295},
      YEAR = {2022},
    NUMBER = {2},
     PAGES = {363--376},
      ISSN = {0025-584X,1522-2616},
   MRCLASS = {46E36 (35K25 42B35 47D07 60J35)},
  MRNUMBER = {4394087},
MRREVIEWER = {Joaquim\ Mart\'{\i}n},
       DOI = {10.1002/mana.202000061},
       URL = {https://doi.org/10.1002/mana.202000061},
}

@book {Ja01,
	AUTHOR = {Jacob, N.},
	TITLE = {Pseudo differential operators and {M}arkov processes. {V}ol.
	{I}},
	NOTE = {Fourier analysis and semigroups},
	PUBLISHER = {Imperial College Press, London},
	YEAR = {2001},
	PAGES = {xxii+493},
	ISBN = {1-86094-293-8},
	MRCLASS = {47G30 (35S05 47D07 60J25)},
	MRNUMBER = {1873235},
	MRREVIEWER = {J. A. van Casteren},
	DOI = {10.1142/9781860949746},
	URL = {https://doi.org/10.1142/9781860949746},
}

\end{document}